\documentclass{article}
\usepackage{amsthm}
\usepackage{amsfonts}

\begin{document}

\centerline{\bf Applications of closed models defined by counting to graph theory and topology.}

\bigskip
\bigskip

{  Tsemo Aristide. }

\bigskip
College Boreal

1, Yonge street, M5E 1E5 Toronto, ON 

Canada

\bigskip
\bigskip

 {\bf Abstract.}  In this paper, we define the notion of  closed models defined by counting, and we compute their homotopy categories.
We apply this construction to various categories of graphs. We show that there does not exist a closed model  in the category
of undirected graphs which characterizes the Ihara Zeta function in the sense that, a morphism $f:X\rightarrow Y$ is a weak equivalence for this model if and only if it induces a bijection between the sets of non degenerated cycles of $X$ and $Y$. Finally, we apply our construction to Galoisian complexes and dessins d'enfant.

\bigskip


\vskip 6mm

\renewcommand{\thefootnote}{}
\footnotetext{ $^*$Corresponding author Tsemo Aristide
\par
E-mail address: tsemo58@yahoo.ca
\par

}

\vskip 6mm
\noindent {\bf\large 1. Introduction  }
\vskip 6mm

The theory of closed models  defined by Quillen  in the context of  category theory provides  foundations of homotopy theory and  is applied to various mathematics areas.
To teach this important idea, it is necessary to have in hands examples which are easy to understand, are not trivial, and can be presented after a short introduction. A good framework to find such closed models is graph theory. In our papers [1] and [2] written in collaboration with Bisson, we have defined such  closed models in the category of directed graphs; one has the virtue to characterize the Zeta function and the other is adapted to symbolic dynamic.
 It is an interesting question to ask whether such a similar closed model  characterizing the Ihara Zeta function
exists in the category of undirected graphs. In this paper, we answer negatively to this question. The closed models defined in [1] and [2]  are particular examples of  closed models defined by counting a family of objects $(X_i)_{i\in I}$ in a topos $C$; that is, the class of weak equivalences of these models is the subclass $W$ of the class of morphisms of $C$, such that for every $f:X\rightarrow Y \in W$,  and every $i\in I$, the  map $c_f:Hom_C(X_i,X)\rightarrow Hom_C(X_i,Y)$ defined by $c_f(h)= f\circ h$ is bijective.
 We start this paper by presenting properties of closed models defined by counting, in particular, we determine their
homotopy categories. We define closed models by counting in subcategories of the category of undirected graphs, which characterize the Ihara Zeta function of  objects of a large subclass of their class of objects. A particular interesting example amongst these closed models is defined in the category $BC_n$, whose objects are $n$-colored graphs.
This category is equivalent to the category of $G_{n}$-sets where $G_{n}$ is the group generated by $a_0,...,a_{n-1}$ such that $a_i^2=1,i=0,...,n-1$.
This category is studied by many  authors, we can quote for example Ladegaillerie [15] who has established an equivalence between $BC_{n+1}$ and the category of  Galoisian $n$-complexes, we deduce the existence of closed models in these categories, and in particular in the category of Galoisian $2$-complexes which is related to dessins d'enfants.

\vskip 6mm
\noindent {\bf\large 1. Plan. }
\vskip 6mm

1. Introduction.

2. Closed models.

3. Closed models defined by counting.

4. Closed models defined by counting in the categories of directed graphs and undirected graphs.

5. Closed models defined by counting in the category of undirected colored graphs.

6. Closed models defined by counting and topology.

\vskip 6mm
{\bf \Large 2. Closed models.}
\vskip 6mm

In this section, we are going to present the basic properties of closed models defined by counting; we start by the following definitions:

\noindent{\bf  Definitions 2.1.}
 A class  $W$  of morphisms of a category $C$ satisfies the $2$-$3$ property if and only if for every
morphisms $f:X\rightarrow Y$  and $g:Y\rightarrow Z$ of $C$,  if two morphisms of the triple $(f,g,g\circ f )$ is an element of $W$, then the third
 is also an element of $W$.

We say that the morphism $g:Y\rightarrow T$ has the right lifting property with respect to $h:X\rightarrow Z$, and that $h$ has the
left lifting property with respect to $g$ if and only if for every commutative diagram:

$$
\matrix{X & {\buildrel{p}\over{\longrightarrow}} & Y\cr h\downarrow && g\downarrow\cr Z
 &{\buildrel{q}\over{\longrightarrow}}& T}
$$

 there exists a morphism $l:Z\rightarrow Y$ such that $l\circ h =p$ and $g\circ l =q$.

 Let $L$ and $R$ be two classes of morphisms of $C$, we say that $(L,R)$ is a weak factorization system if and only if:

 - Every morphism $f\in C$ can be written $f = r\circ l$ where $l\in L$, $r\in R$;

 - $L$ is the class of morphisms which has the left lifting property in respect of every morphism of $R$;

 - $R$ is the class of morphisms which has the right lifting property in respect of every morphism of $L$.

\medskip

Let $C$ be a category complete and cocomplete; we say that $C$ is endowed with a closed  model  if and only if there exist three classes of morphisms $(Fib,Cof,W)$ such that:

- $W$ satisfies the $2$-$3$ property,

- Let $Fib'=W\cap Fib$,  $(Cof,Fib')$ is a weak factorization system

- Let $Cof' = W\cap Cof$, $(Cof',Fib)$ is a weak  factorization system.

\vskip 2mm

We start by the following general example:

\medskip

\noindent{\bf Proposition 2.2.}
{\it Let $C$ be a category complete and cocomplete, let $W$ be a class of morphisms of $C$ which satisfies the $2$-$3$ property.
Suppose that there exists a class of morphisms $Cof$ of $C$ such that $(Cof,W)$ is a weak factorization system. Then, there exists
a closed model on $C$ whose class of weak equivalences is $W$, its class of cofibrations is $Cof$, its class of weak cofibrations $Cof'$ is the
 class $Iso(C)$ of isomorphisms of $C$,  its class of fibrations $Fib$ is
the class $Hom(C)$ of all morphisms of $C$, and its class of weak fibrations $Fib'$ is $W$.}

\vskip 2mm
\noindent{\bf Proof.}
We have:
 $(Cof',Fib) = (Iso(C),Hom(C))$ and $(Cof,Fib')=(Cof,W)$ are weak factorization systems.  We also have $Fib\cap W = Hom(C)\cap W =W =Fib'$,
 and $Cof\cap W = Iso(C)$ Since $(Cof,W)$ is a factorization system.
  We deduce  that $(Hom(C),Cof,W)$ defines a closed  model on $C$.
\vskip 2mm 

\vskip 6mm
{\bf \Large 3. Closed model defined by counting.}
\vskip 6mm

We are going to apply the previous result to define  closed models to count objects in categories.
Let $C$ be a category complete and cocomplete whose initial object is denoted by $\phi$.
 For every objects $X$ and $Y$ of $C$, we denote by $X+Y$ the sum of $X$ and $Y$.
 Let $(X_i)_{i\in I}$ be a family of objects of
$C$ and $l_i: \phi\rightarrow X_i$ the canonical morphism.
There exist morphisms $j^i_1:X_i\rightarrow X_i+X_i$ and $j^i_2:X_i\rightarrow X_i+X_i$ such that for every
morphisms $f:X_i\rightarrow Z$ and $g:X_i\rightarrow Z$, there exists a unique morphism $m(f,g):X_i+X_i\rightarrow Z$
such that $m(f,g)\circ j^i_1 = f$ and $m(f,g)\circ j^i_2 = g$. We set $m_i = m(Id_{X_i},Id_{X_i})$.
Such a morphism is often called a folding morphism. We suppose that the class of morphisms $l_i,m_i\in I$ admits the small element argument (see [11] 12.4.13).
We denote by $W_I$ the class of morphisms which are right orthogonal to every  morphisms $l_i$ and $m_i, i\in I$.

\medskip

\noindent{\bf Proposition 3.1.}
{\it  A morphism $f:X\rightarrow Y $ of $C$ is an element of $W_I$ if and only if for every $i\in I$, the map $c^i_f:Hom_C(X_i,X)\rightarrow Hom_C(X_i,Y)$
defined by $c^i_f(h) = f\circ h$ is bijective. We deduce that  $W_I$ satisfies the $2$-$3$-property and there exists a closed model on $C$ whose class of weak equivalences is $W_I$.}

\vskip 2mm
\noindent{\bf Proof.}
Let $f:X\rightarrow Y$ be a morphism of $C$, suppose that for every $i\in I$, $f$ is orthogonal to $l_i$ and $m_i$.
Let $h,h'\in Hom_C(X_i,X)$ such that $f\circ h = f\circ h'$.  The following diagram commutes:

$$
\matrix{X_i+X_i &{\buildrel{h+h'}\over{\longrightarrow}}& X\cr m_i\downarrow  && f\downarrow\cr
X_i & {\buildrel{f\circ h}\over{\longrightarrow}} & Y}
$$

Since $f$ is right orthogonal to $m_i$, we deduce the existence of a morphism $l:X_i\rightarrow X$ such that
$l\circ m_i =h+h'$.
We have $l\circ m_i\circ j^i_1 =l\circ m_i\circ j^i_2 =l$. We deduce that $ h =(h+h')\circ j^i_1 = l\circ m_i\circ j^i_1
= l\circ m_i\circ j^i_2= (h+h')\circ j^i_2 = h'$.

Let $h:X_i\rightarrow Y$ be any morphism, the following diagram commutes:

$$
\matrix{ \phi & \longrightarrow & X\cr \downarrow && \downarrow f\cr
X_i & {\buildrel{h}\over{\longrightarrow}} & Y}
$$
thus it has a filler $p:X_i\rightarrow X$ such that $f\circ p = h$.
This implies that $c_f^i$ is bijective.

We show now that $W_I$ satisfies the $2$-$3$ property. Let $f:X\rightarrow Y$ and $g:Y\rightarrow Z$ be
morphisms of $C$, $c^i_{g\circ f}= c^i_g\circ c^i_f$. Since $c^i_f,c^i_g$ and $c^i_{g\circ f}$ are morphisms of sets,
we deduce that if two morphisms of the triple $(c^i_f,c^i_g,c^i_{g\circ f})$ are bijective, so is the third.

Let $cell(I)$ be the class of morphisms of $C$ which are  retracts of  transfinite compositions of pushouts of $l_i,m_i,i\in I$,  the propositions
 12.4.14 and 12.4.20 of [11] imply that $(cell(I),W_I)$ is a factorization system. We deduce from the proposition 2.2
the existence of a closed model on $C$, whose class of weak equivalences is $W_I$.
\vskip 2mm

\medskip

\noindent{\bf Remarks.}
Recall that  a closed model $(Fib,Cof,W)$ on a category is cofibrantly
generated (see [11] 13.2.2.) if there exist  sets of morphisms $(f_i)_{i\in I}$ and $(g_j)_{j\in J}$ both which allow the small element argument, such that
the class of fibrations is the class of morphisms which are right orthogonal to every morphism of the family $(g_j)_{j\in J}$, and the class of weak fibrations is the class of morphisms which are
right orthogonal to every morphism of the family $(f_i)_{i\in I}$. Thus, the closed model defined by counting in the previous proposition is cofibrantly generated. In the sequel, we will only consider such closed models defined by counting.

Let $K$ be a subset of $I$, we denote by $X_K$ the sum of objects of the family $(X_k)_{k\in K}$. The morphism $\phi\rightarrow X_K$ is a cofibration,
since it is transfinite composition of pushouts of elements of $(l_k)_{k\in K}$.

\vskip 2mm
\noindent {\bf The homotopy category of a closed model defined by counting.}
\vskip 2mm

 One of the main purpose of the theory of closed models is to find a proper framework to localize classes of morphisms.
 In this perspective, we are going to compute the homotopy category of a closed model defined by counting. We start by remarking
the fact that since every morphism is a fibration in a closed model defined by counting,   every object is fibrant. Let us determine now
the  cofibrant replacement of an object:

\medskip

\noindent{\bf Proposition 3.2.}
{\it Let $U$ be a Grothendieck universe, and $C$ be an $U$-category endowed with a closed model defined by counting  a set of objects $(X_i)_{i\in I}$, where $I$ is an $U$-set.
 Let $Z$ be an object of  $C$, then $Z$ has a cofibrant
replacement $QZ$, isomorphic to a transfinite composition of a subset of $l_i,m_i,i\in I$.}

\vskip 2mm
\noindent{\bf Proof.}
Let $I_Z=\{f\in Hom_C(X_i,Z),i\in I\}$. It is an $U$-set.  We  denote by $X_{I_Z}$
 the pushout of elements of $I_Z$, and by $i_f:X_i\rightarrow X_{I_Z}$ the morphisms satisfying the universal property of the pushout. There exists a morphism $d_Z:X_{I_Z}\rightarrow Z$ such that for every $f\in I_Z$, $d_Z\circ i_f = f$.
   For every $i\in I$, $c_{d_Z}:Hom_C(X_i,X_{I_Z})\rightarrow Hom_C(X_i,Z)$ defined by $c_{d_Z}(g)=d_Z\circ g$
 is surjective since if  $h\in Hom_C(X_i,Z)$, we have $d_Z\circ i_h =h$. Consider  $L(Z)$ the $U$-set whose elements are morphisms $p_V:X_{I_Z}\rightarrow V$
 such that $p_V$ is a transfinite composition of pushouts of a subset of  $m_i, i\in I$, and there exists a morphism
 $f_V:V\rightarrow Z$ such that $f_V\circ p_V = d_Z$. There exist a relation of order define on $L(Z)$ such that $p_V\geq p_W$
  if and only if there exists a morphism $h_{V,W}:W\rightarrow V$ such that $p_V=h_{V,W}\circ p_W$. Let $(p_{V_j})_{j\in J}$
  be an ordered family of $L(Z)$; $lim_{j\in J}p_{V_j}$ is a lower bound of $(p_{V_j})_{j\in J}$. The Zorn's lemma implies  that
    $L(Z)$ has a maximal $c_Z:QZ\rightarrow Z$ which is a weak equivalence.
\vskip 2mm

\noindent{\bf Remark.}
In the rest of this section, we are going to suppose that the category $C$ is $U$-small, where $U$ is a Grothendieck universe. Let $V$ and $W$ be objects of $C$,  every morphism $f:V\rightarrow W$ induces a morphism $d(f):X_{I_V}\rightarrow X_{I_W}$;
 $d(f)$ also induces a morphism $Qf:QV\rightarrow QW$ such that $f$ is a weak equivalence
 if and only if $Qf$ is an isomorphism. (See also [6] Lemma 5.1).

\medskip

 \noindent{\bf Definitions 3.3.}
A path object of $Z$ is an object $Z^I$ such that there exists a weak equivalence $i_Z:Z\rightarrow Z^I$, a morphism $p_Z:Z^I\rightarrow Z\times Z$
such that $(id_Z,id_Z) = p_Z\circ i_Z$.

Two morphisms $f,g:Y\rightarrow Z$ are right homotopic if and only if there exists a path object $Z^I$, and a morphism $H:Y\rightarrow Z^I$
such that $(f,g) = p_Z\circ H$.

\medskip

\noindent{\bf Proposition 3.4.}
{\it Let $C$ be a category endowed with a closed model defined by counting, and  $Y$ a cofibrant
object of $C$. Two morphisms $f,g:Y\rightarrow Z$ of $C$ are right homotopic if and only if they are equal.}

\vskip 2mm
\noindent{\bf Proof.}
If $Y$ is cofibrant and $f,g$ are right homotopic, we can suppose that there exists a path object $Z^I$ such that
$i_Z:Z\rightarrow Z^I$ is an acyclic cofibration and $(f,g) = p_Z\circ H$ (See [6] Lemma 4.15).
 Since $i_Z$ is an acyclic cofibration of a closed model defined by counting, it is
an isomorphism. We can thus suppose that $Z^I = Z$ and $p_Z=(id_Z,id_Z)$. This implies that $f=p_1\circ p_Z\circ H = p_2\circ p_Z\circ H =g$
where $p_1,p_2:Z\times Z\rightarrow Z$ are the projections on the first and second factors.
\vskip 2mm

\noindent{\bf Remark.}
Let $Y$ and $Z$ be two objects of $C$, we denote by $C(Y)$ (resp. $C(Z)$) the cofibrant replacement of $Y$ (resp. $Z$). The objects
$C(Y)$ and $C(Z)$ are also fibrant. The homotopy category of $C$ is the category which have the same class of objects than $C$,
and the set $Hom_{Hot}(Y,Z)$ of morphisms of the homotopy category between the objects $Y$ and $Z$
is the set $\pi(X,Y)$ whose elements are  right homotopy classes of morphisms between $C(Y)$ and $C(Z)$. (See [6] 4.22 and definition 5.6).
We deduce:

\medskip

\noindent{\bf Proposition 3.5.}
{\it Let $C$ be a category endowed with a closed model defined by counting, for every objects $X$ and
$Y$ of $C$, we have $Hom_{Hot}(Y,Z)= \pi(C(Y),C(Z)) = Hom_C(C(Y),C(Z))$.}

\vskip 6mm
{\bf \Large 4.  Closed models defined by counting in the categories of directed graphs and undirected graphs.}
\vskip 6mm

We will define now various closed models in different categories of graphs. We start by the category of
directed graphs.

Let $C_D$ be the category which has two objects that we denote by $0$ and $1$. We suppose
that $Hom_{C_D}(0,1)$ contains two elements $s,t$, $Hom_{C_D}(0,0)$,  $Hom_{C_D}(1,1)$ contain one element and $Hom_{C_D}(1,0)$ empty.

\medskip

\noindent{\bf Definition 4.1.}
A directed graph is a presheaf defined on $C_D$. Let $X$ be such a presheaf; $X$ is defined by two sets
$X(0)$ and $X(1)$, and two maps $X(s),X(t):X(1)\rightarrow X(0)$.

The set $X(0)$ is called the space of nodes of $X$ and the set $X(1)$ the space of directed arcs of $X$.

\medskip

\noindent{\bf Definition 4.2.}
A morphism $f:X\rightarrow Y$ between the graphs $X$ and $Y$ is a natural transformation between the
presheaves $X$ and $Y$; thus $f$ is defined by a morphism $f_0:X(0)\rightarrow Y(0)$, $f_1:X(1)\rightarrow Y(1)$
such that $f_0 \circ X(s) = Y(s)\circ f_1$, $f_0\circ X(t) = Y(t)\circ f_1$ .

We denote by $Gph$ the category of directed graphs.
The category of  directed graphs is complete and cocomplete since it is a Grothendieck topos.

\vskip 2mm

Examples of directed graphs are: the directed dot graph $D$. It is the graph such that $D(0)$ is a singleton
and $D(1)$ is empty.

The directed arc graph $A$ is the graph defined by: $A(0)=\{u,v\}$, $A(1)=\{a\}$ and $A(s)(a)=u$, $A(t)(a)=v$.

Let $p$ be a strictly positive integer. We denote by $c_p$ the graph whose set of nodes is  $Z/pZ$, let $[n]$ be the class of the integer
$n$ in $Z/p$, there exists a unique arc $a_n$ such that $c_p(s)(a_n) =[n]$ and $c_p(t)(a_n) =[n+1]$.

We can define on $Gph$ the closed model obtained by counting the elements of the set $Cycl=\{c_p,p\in N-\{0\}\}$.

\medskip

\noindent{\bf Remark.}
The closed model obtained here have the same class of weak equivalences than the closed model defined in [1], but  the classes
of cofibrations, weak cofibrations, fibrations, weak fibrations of these closed models are different.

\medskip
\noindent{\bf Definition 4.3.}
Let $X$ be a  directed graph, for every non zero integer $p$,  we denote by $n_p(X)$ the cardinality of $Hom_{Gph}(c_p,X)$.
Suppose that for every strictly positive integer $p$, $n_p(X)$ is finite. The Zeta serie $Z_X(t)$ of $X$ is:

$$
exp(\sum_{p=1}^{p=\infty}n_p(X){{t^p}\over p}).
$$

\medskip

\noindent{\bf Proposition 4.4.}
{\it Let $X$ and $Y$ be two finite directed graphs, $Z_X(t) = Z_Y(t)$ if and only if there exists an isomorphism
$f$ in $Hom_{Hot}(X,Y)$.}

\vskip 2mm

\noindent{\bf Proof.}
Let $X$ and $Y$ be two finite graphs; the proposition 3.5 implies that there exists an isomorphism in $Hom_{Hot}(X,Y)$ if and only if there
exists an isomorphism of graphs $f:C(X)\rightarrow C(Y)$. The proposition 3.2 implies that the cofibrant replacement $C(X)$ is a sum of cycles such that $Z_X(t)=Z_{C(X)}(t)$. We deduce that $Z_X(t)=Z_{C(X)}(t) = Z_{C(Y)}(t) = Z_Y(t)$.

\vskip 2mm
{\bf Undirected graphs.}
\vskip 2mm

It is natural to try to generalize this closed model to others categories of graphs, unfortunately straightforward
generalizations do not have the same natural properties, for example we do not obtain the same characterization
of the weak equivalences with the corresponding Zeta series. We will  consider the category $UGph$ of undirected graphs.

Let $C_U$ be the category which has two objects that we denote by $0$ and $1$. We suppose
that $Hom_{C_U}(0,1)$ contains two elements $s,t$, $Hom_{C_U}(0,0)$ contains one element,
 $Hom_{C_U}(1,1)$ contains the identity and an involution $i$ such that $i\circ s =t$, and $Hom_{C_U}(1,0)$ is empty.

\medskip

\noindent{\bf Definition 4.5.}
An undirected graph is a presheaf defined on $C_U$. Let $X$ be such a presheaf, $X$ is defined by two sets
$X(0)$ and $X(1)$, two maps $X(s),X(t):X(1)\rightarrow X(0)$ and an involution $X(i)$ of $X(1)$ such that $X(s)\circ X(i) = X(t)$.

The set $X(0)$ is called the space of nodes, and the space $X(1)$ the space of half-arcs. For an half-arc $a\in X(1)$,
$X(s)(a)$ is the source of $a$ and $X(t)(a)$ is the target of $a$.

Remark that $X(i)$ is an involution of $X(1)$,
and the source of the half-arc $a$ is the target of $X(i)(a)$ since $i\circ s =t$.

We have not assume that $X(i)$ acts freely, this implies the existence of undirected graphs $X$ with degenerated loops; these are
half-arcs fixed by $X(i)$.

An arc of the graph $X$ is defined by a couple $(u,X(i)(u))$ where $u\in X(1)$. We denote by $Arc(X)$ the space
of arcs of the undirected graph $X$. The source or the target of $u$ will often be called an end of the arc $(u,X(i)(u))$.

Geometrically, if the set of half arcs of an undirected graph $X$ does not contain  a degenerated loop, it can be represented by a set of points corresponding to its nodes,
and an arc $(u,X(i)(u))$ is an unoriented segment connecting $X(s)(u)$ and $X(t)(u)$.

\medskip

\noindent{\bf Definition 4.6.}
A morphism $f:X\rightarrow Y$ between the undirected graphs $X$ and $Y$ is a natural transformation between the
presheaves $X$ and $Y$; thus $f$ is defined by  morphisms $f_0:X(0)\rightarrow Y(0)$, and $f_1:X(1)\rightarrow Y(1)$
such that $f_0 \circ X(s) = Y(s)\circ f_1$, $f_0\circ X(t) = Y(t)\circ f_1$ and $f_1\circ X(i) = Y(i)\circ f_1$.

The morphism of graphs $f:X\rightarrow Y$ induces a morphism $a(f):Arc(X)\rightarrow Arc(Y)$. If there is no confusion, we will often
denote $a(f)$ by $f_1$.

Examples of undirected graphs are: the undirected dot graph $D_U$. It is the graph such that $D_U(0)$ is a singleton
and $D_U(1)$ is empty.

The undirected arc graph $A_U$ is the graph defined by $A_U(0)=\{u_1,u_2\}$, $A_U(1)=\{a_1,a_2\}$ such that $A_U(i)(a_1)=a_2$, $A_U(s)(a_1)=u_1$
and $A_U(t)(a_1)=u_2$.

The graph $V_U$ is the graph defined by $V_U(0)=\{v_1,v_2,v_3\}$, $V_U(1)=\{b_1,b_2,c_1,c_2\}$ such that $V_U(s)(b_1)=V_U(s)(c_1)=v_1$,
$V_U(t)(b_1) = v_2, V_U(t)(c_1)=v_3$, $V_U(i)(b_1)=b_2$ and $V_U(i)(c_1) = c_2$.

The path graph $P_n$ is the graph whose set of nodes is $\{0,...,n\}$, $P_n(1)=\{p^+,p^-, p= 0,...,n-1\}$ such that
$P_n(s)(p^+) = p, P_n(t)(p^+)=p+1, P_n(i)(p^+) = p^-$.

There is a  morphism $f:V_U\rightarrow A_U$ such that $f_0(v_1)=u_1, f_0(v_2)=f_0(v_3)= u_2$, $f_1(b_1)=f_1(c_1) = a_1$.
This morphism is called the elementary folding.

Let $p$ be a strictly positive integer. We denote by $c^p_U$ the undirect graph whose set of nodes is $Z/pZ$. Let $[n]$ be the class of the integer
$n$ in $Z/pZ$, we have $c^p_U(1)=\{[n]^+,[n]^-, [n]\in Z/p\}$, $c^p_U(s)([n]^+)=[n], c^p_U(s)([n]^-)=[n+1]$, and $c^p_U(i)([n]^+)=[n]^-$.  The graph $c^p_U$ is called the undirected $p$-cycle.

Let $X$ and $Y$ two undirected graphs isomorphic to the $1$-cycle $c_U^1$. Remark that there exists a unique morphism  $f:D_U\rightarrow X$ 
(resp. $g:D_U\rightarrow Y$).  The pushout of $f$ and $g$ is
the eight graph. Geometrically, it corresponds to two circles attached in one point.

\medskip

\noindent{\bf Definition. 4.7.}
Let $X$ be an undirected graph, a $p$-cycle of $X$ is a morphism $f:c^p_U\rightarrow X$. We say that the $p$-cycle
$f$  has a backtracking if and only if there exists an integer $n$ such that $f_1([n+1]^+)=f_1([n]^-)$. We denote by $Cycl_p(X)$
the set of $p$-cycles of the undirected $X$ without a backtracking.

We denote by $W_U$ the class of morphisms of $UGph$ such that for every $f:X\rightarrow Y$ in $W_U$, for every integer $p>0$, the morphism
$c_p(f):Hom(c^p_U,X)\rightarrow Hom(c^p_U,Y)$ which sends the morphism $h$ to $f\circ h$ induces a bijection on cycles without a backtracking.
The class of morphisms $W_U$ satisfies the $2$-$3$-property.

Let $X$ be a finite undirected graph, we denote by $c_p(X)$ the cardinality of the set of morphisms $c^p_U\rightarrow X$
without a backtracking. The Ihara zeta function of $X$ is defined by:

$$
exp(\sum_{p\geq 1}{{c_p(X)}\over p}t^p)
$$

Remark that if $f:X\rightarrow Y$ is a morphism between finite undirected graphs in $W_U$,  the graphs $X$ and $Y$ have the same Ihara zeta function.
We want to find a closed model for which $W_U$ is the class of weak equivalences. We will see that such a model does not exist.

\medskip

\noindent{\bf Remark.}
We can naively adapt the previous closed model defined in the category of directed graphs to the category of undirected graphs: so, we define the closed model
obtained by counting elements of the family $(c^p_U)_{p\in N-\{0\}}$. A morphism $f:X\rightarrow Y$ of $UGph$
is a weak equivalence for this closed model if and only if  for every strictly positive integer $p$, the map $Hom(c^p_U,X)\rightarrow Hom(c^p_U,Y)$ induced by $f$
is bijective. Thus, $f$ induces a bijection between the $p$-cycles of $X$, and the $p$-cycles of $Y$ for every strictly
positive integer, but  the image of a cycle without backtracking by $f$ is not necessarily a cycle without backtracking. This closed model is essentially trivial as shows the following result:

\medskip

\noindent{\bf Proposition 4.8.}
{\it A weak equivalence $f:X\rightarrow Y$ between two undirected finite and connected graphs for the closed model defined
by counting elements of the family $(c^p_U)_{p\in N-\{0\}}$ is an isomorphism.}

\vskip 2mm

{\bf Proof.}
Firstly, we show that $f_1:X(1)\rightarrow Y(1)$ is injective. Let $a,b$ be two distinct arcs such that $f_1(a)=f_1(b)$.
 There exist morphisms $h_i:c^2_U\rightarrow X$, $i=1,2$ such that the image of $h_1$ is $a$ and the image of $h_2$ is $b$,
 the square diagram:

 $$
 \matrix{ c^2_U+c^2_U & {\buildrel{h_1+h_2}\over{\longrightarrow}}& X\cr
 \downarrow j_2 &&\downarrow f\cr
 c^2_U &{\buildrel{f\circ h_1}\over{\longrightarrow}} & Y}
 $$

 is commutative and does not have a filler this is a contradiction.

We show now that $f_1$ is surjective. Let $a$ be an arc of $Y$, there exists a morphism $h:c_U^2\rightarrow Y$ whose image is $a$.
Since $f$ is right orthogonal to $i_2:\phi\rightarrow c^2_U$, we deduce the existence of a morphism $h':c^2_U\rightarrow X$ such that
$h = f\circ h'$. This implies that $f_1$ is surjective on arcs.

We show now that $f_0:X(0)\rightarrow Y(0)$ is injective. Let $x$ and $y$ be two distinct nodes of $X$
such that $f_0(x)=f_0(y)$. Since $X$ is connected, there exists a path $h$ between $x$ and $y$, $f(h)$ is a $p$-cycle where $p>0$.
Since $f$ induces a bijection on $p$-cycles, we deduce the existence of a $p$-cycle $c$ of $X$ such that $f(c) =f(h)$. This is in contradiction
with the fact that $f$ is bijective on arcs.

We show now that $f_0$ is surjective.

Let $y$ be a node of $Y$ since $Y$ is connected, there exists an arc $b$ of $Y$ which has $y$ as an end. Since $f_1$ is
bijective, we deduce the existence of an arc $a$ of $X$ such that $b = f_1(a)$. This implies that $a$ has an end $x$ such that $f_0(x)= y $.
\vskip 2mm

\noindent{\bf Theorem 4.9.}
{\it There does not exist a closed model on $UGph$ whose class of weak equivalences is the class $W_U$.}

\vskip 2mm

\noindent{\bf Proof.}
Suppose that such a closed model exists, then the elementary folding $f:V_U\rightarrow A_U$ would be a weak
equivalence, and we can write $f=g\circ h$ where $h:V_U\rightarrow X$ is a weak cofibration and $g: X\rightarrow A_U$ is a fibration. The $2$-$3$ property
implies  that $g$ is a weak fibration.

Suppose that $h_1(b_1) = h_1(c_1)$, and consider the morphism $l:V_U\rightarrow c^2_U$ defined by $l_0(v_1) = [0], l_0(v_2)=l_0(v_3)=[1]$,
$l_1(b_1)= [0]^+$ and $l_1(c_1)=[1]^-$.
 We can define the pushout diagram:

$$
\matrix{V_U & {\buildrel{h}\over{\longrightarrow}} & X\cr l\downarrow && p\downarrow\cr c^2_U
 &{\buildrel{q}\over{\longrightarrow}}& Z}
$$

 Remark that $X$ does not have any cycle without backtracking since $h$ is a weak equivalence
and $V_U$ does not have any cycle without backtrackings. We deduce that $Z$  does not have any cycle without backtrackings since $Z$ is isomorphic to $X$. This implies that $q$ is not a weak equivalence.
This is a contradiction with the fact that in a closed model, the pushout of a weak cofibration is a weak cofibration.
Thus $h_1(b_1)$ is distinct of $h_1(c_1)$. Remark that the image of $V_U$ by $h$ cannot be isomorphic to a $2$-cycle since $h$ is a
weak equivalence; we deduce that this image is isomorphic to $V_U$.

Now, consider the morphism $m:c^2_U\rightarrow A_U$ defined by $m_0([0]) = u_1, m_0([1])=u_2$ and $m_1([0]^+)=m_1([1]^-)=a_1$. Consider the pullback diagram:

$$
\matrix{U & {\buildrel{p'}\over{\longrightarrow}} & X\cr q'\downarrow && g\downarrow\cr c^2_U
 &{\buildrel{m}\over{\longrightarrow}}& A_U}
$$

The morphism $q':U\rightarrow c^2_U$ must be a weak equivalence since in a closed model, the pullback of a weak fibration
must be a weak fibration.  But, there exists a subgraph of $U$ isomorphic to
the pullback of the elementary folding by  $m:c^2_U\rightarrow A_U$. Such a subgraph   is
isomorphic to the  graph obtained by identifying two nodes of two distinct unoriented $2$-cycles. So $q'$ cannot be a weak equivalence. This is a contradiction.

\vskip 6mm
{\bf \Large 5. Closed models defined by counting in the category of undirected colored graphs.}
\vskip 6mm

We are going to define  closed models in  subcategories of $UGph$, and in particular in the category
of undirected colored graphs.

\vskip 6mm

{\bf Definitions 5.1.}
Let $X$ be an undirected graph, for every node $x$ of $X$, we denote by $X(x,*)$ the set of arcs $(u,X(i)(u))$
such that $X(s)(u) =x$ or $X(t)(u)= x$.

A morphism $f:X\rightarrow Y$ between undirected graphs is  a covering if and only if for every $x\in X_0$, the morphism
$f_x:X(x,*)\rightarrow Y(f_0(x),*)$ induced by $f$ is bijective. Remark that for every undirected graph, the morphism $\phi\rightarrow X$
is a covering.

Let $X$ be an undirected graph, we denote by $C_X$ the category whose objects are coverings $f:Y\rightarrow X$.
A morphism between the objects $f:Y\rightarrow X$ and $g:Z\rightarrow X$ is a covering morphism $h:Y\rightarrow Z$ such
that $g\circ h = f$.

\medskip

\noindent{\bf Proposition 5.2.}
{\it  Limits and colimits exist in $C_X$.}

\vskip 2mm

\noindent{\bf Proof.}
The class $Cov$ of covering morphisms of $UGph$ is the class of morphisms which are right orthogonal to the elementary
folding $f_U:V_U\rightarrow A_U$ and to $i_U:D_U\rightarrow A_U$. This implies that the pullback of a covering morphism is a covering morphism. Since
the products in $C_X$ are pullbacks of covering morphisms, we deduce that products and pullbacks exist in $C_X$, and henceforth that
limits exist in $C_X$. (See SGA 4.1 proposition 2.3).

Let $(f_i:Y_i\rightarrow X)_{i\in I}$ be a family of elements of $C_X$. The morphism $f:\sum_{i\in I}Y_i\rightarrow X$
whose restriction to $Y_i$ is $f_i$ is a covering; this implies that sums exist in $C_X$.

We show now that pushouts exist in $C_X$. Let $h:Z\rightarrow X$ and $h':Z'\rightarrow X$ be two objects of $C_X$;
consider an object $p:Y\rightarrow X$,  $f:Y\rightarrow Z$ and $g:Y\rightarrow Z'$ two morphisms of $C_X$.
Without restricting the generality, we can suppose that $Y,Z$ and $Z'$ are connected. The morphisms $f$ and $g$ are surjective
on nodes and arcs since they are coverings and the pushout of $f$ and $g$ is defined by the graph $L$  whose set of nodes,
$L_0$ is the quotient of $Z_0\bigcup Z'_0$ by the equivalence relation generated by: let $x\in Z_0$ and $x'\in Z'_0$, $x\simeq x'$ if and only if there exists $x"\in Y_0$ such that
  $f_0(x") = x$ and $g_0(x") =x'$;  $L_1$ is the quotient of $Z_1\bigcup Z_1'$ by the equivalence relation generated by: let $a\in Z_1$ and $b\in Z'_1$,
$a\simeq b$ if and only if there exists $c\in Y_1$ such that $f_1(c) =a$ and $g_1(c) =b$. We denote by $p_l:Y\rightarrow L$ the quotient morphism. There exists a  morphism $l:L\rightarrow X$ such that $p = l\circ p_l$
since $ p = h\circ f = h'\circ g$. The morphism $l$ is a covering since $f$ and $g$ are elements of $C_X$; it  is the pushout of $f$ and $g$.

\vskip 2mm

We consider $R_X^p$ the set of graphs in $C_X$ such that an element of $R_X^p$ is obtained by attaching a forest
to a $p$-cycle. We define on $C_X$ the closed model obtained by counting elements of $R_X=\{R_X^p,p\in N-\{0\}\}$.
Thus a morphism $f:X\rightarrow Y$ of $C_X$ is a weak equivalence if and only if for every for every $p>0$, for every $U_p\in R^p_X$,
$Hom_{C_X}(U_p,X)\rightarrow Hom_{C_X}(U_p,Y)$ is bijective.

\medskip

\noindent
{\bf Proposition 5.3.}
{\it Let $f:Y\rightarrow X$ be an object of $C_X$ such that $Y$ does not have loops, then for every $p$-cycle  $h:c^p_U\rightarrow Y,p>1$,
without backtrackings there exists an element $U$ of $R_X^p$ and a morphism $g:U\rightarrow Y$ of $C_X$ whose restriction to the $p$-cycle is $h$.}

\vskip 2mm

\noindent{\bf Proof.}
Let $h:c^p_U\rightarrow Y$ be a $p$-cycle. Since $Y$ does not have loops, for every integer $n$, $h_1([n]^+)$ is distinct
of $h_1([n+1]^+)$. This implies that we can attach a tree to every node of $c_p^U$ to obtain a graph $U$ for which there exists a
covering $g:U\rightarrow Y$ whose restriction to $c_p^U$ is $h$.

\vskip 2mm

{\bf Remark.}
The previous proposition is not true if there are $1$-cycles in $Y$. Consider the following example: $X=Y$ is the $1$-cycle,  Consider the
 morphism  defined by  $h:c^2_U\rightarrow X$ such that
$h_1([0]^+)=h_1([1]^+)$. This cycle does not have backtracking, but it is impossible extend $h$ to an element of $R_X^2$ such that it becomes a covering,
since the restriction of $h_1$ to $c^2_U([1],*)$ is not injective.

\medskip

\noindent{\bf Proposition 5.4.}
{\it Let $f:Y\rightarrow X$ and $g:Z\rightarrow X$ be objects of $C_X$ such that $Y$ and $Z$ are finite and do not have loops. If there
exists a weak equivalence between $f$ and $g$,
then $Y$ and $Z$ have the same Ihara Zeta function.}

\vskip 2mm

\noindent{\bf Proof.}
A weak equivalence between $f:Y\rightarrow X$ and $g:Z\rightarrow X$ is defined by a covering $h:Y\rightarrow Z$
such that $g\circ h = f$. We are going to show that $h$ induces a bijection on $p$-cycles without backtracking. Let $u,u':c^p_U\rightarrow Y$ be two $p$-cycles
of $Y$ without backtrackings such that $h\circ u = h\circ u'$. Since $Z$ does not have loops, there exists an element $v:V\rightarrow X\in R_X^p$
 and a morphism between $v$ and $g$ whose restriction to the $p$-cycle of $V$ coincide with $h\circ u$. Since $h$ is a covering, we can lift $v$
to morphisms $v_1:V\rightarrow Y$ (resp $v_2:V\rightarrow Y$) whose restriction the $p$-cycle is $u$ (resp. $u'$) and such that $v_1,v_2$ are morphisms
of $C_X$ respectively between $v$ and $f$.
Since $h$ is a weak equivalence, we deduce that $v_1=v_2$, and henceforth that $u= u'$. We deduce that
 $h$ is injective on $p$-cycles. The fact that $f$ is surjective on $p$-cycles results from the fact that
for every $p$-cycle without backtracking $l:c^p_U\rightarrow Z$ there exists an element $v:V\rightarrow X\in R_X^p$ and a morphism $d$ of $C_X$ between
$v$ and $g$ whose restriction to the $p$-cycle of $V$ is $l$. We can lift $d:V\rightarrow Z$ to  a morphism $d':V\rightarrow Y$ since $h$ is a covering, the restriction
of $d'$ to the $p$-cycle of $V$ is a preimage of $l$.

\vskip 2mm

\noindent{\bf Remarks.}
Let $f:Y\rightarrow X$ be an object of $C_X$ without loops. A  $p$-cycle $u:c_U^p\rightarrow Y$ without backtracking, is primitive if and only if for every $q$-cycle $u':c_U^q\rightarrow Y$ without backtracking, if there exists a morphism $f:c_U^p\rightarrow c_U^q$ such that $u=u'\circ f$, then $p=q$. Two primitive $p$-cycles $u,u':c_U^p\rightarrow Y$ are equivalent if there exists an isomorphism $f$ of $c_U^p$ such that $u'=u\circ f$. We denote by $E_p(Y)$ the set whose elements are equivalence classes of primitive $p$-cycles without backtracking, and by $E(Y)=\bigcup_{p\in N-\{0\}}E_p(Y)$. For a primitive $p$-cycle without backtracking $u:c_U^p\rightarrow Y$, we will denote by $[u]$ its equivalence class. 

For every element $[u]\in E_p(Y)$, we choose an element  $u:c_U^p\rightarrow Y$ in this class, and consider the element
$V_u$ of $R^p_X$ such that there exists a morphism $v_u:V_u\rightarrow Y$ of $C_X$ whose restriction to the cycle of $V_u$
coincide with $u$. Let $V$ be the direct summand of the
graphs $V_u$, there exists a covering $c: V\rightarrow Y$ whose restriction to $V_u$ is $v$. The morphism $c$ is
   a weak equivalence and the morphism $\phi\rightarrow c(Y)$ is a
cofibration, thus $c(Y)$ is a cofibrant replacement of $Y$.

\vskip 2mm

Let $B_n$ be the undirected graph which has one node $*$, and $n$ undirected loops. Let $X$ be an undirected graph, there
exists a covering $f:X\rightarrow B_n$ if and only if $X$ is a $n$-regular graph and the edges of $X$ can be colored by $n$-distinct colors.
The proposition 5.4 implies that there exists a closed model on the category of $n$-regular graphs whose edges can be colored
by $n$ distinct colors such that, if there exists a weak equivalence between two finite graphs without loops in this category, then they
have the same Zeta function.

\medskip

\noindent{\bf Proposition 5.5.}
{\it  Suppose that $X=B_n$, let $f:Y\rightarrow X$ and $g:Z\rightarrow X$ be two objects of $C_X$ such that $Y$ and $Z$ are finite
and do not have loops. Moreover, suppose that $Y$ and
$Z$ have the same Zeta function, and there exist an  isomorphism $H_p:E_p(Y)\rightarrow E_p(Z)$  such that each element $[c]\in E_p(Y)$ there exists a morphism $c:c_U^p\rightarrow Y$ representing $[c]$, such that
$f\circ c = g\circ H_p(c)$ where $H_p(c)$ represents $H_p([c])$. Then there exists a graph $L$, coverings $p:L\rightarrow Y$ and $p':L\rightarrow Z$ which
are weak equivalences.}

\vskip 2mm

\noindent{\bf Proof.}
If $Y\rightarrow X$ and $Z\rightarrow X$ are two objects of $C_X$ such that the Zeta series of $Y$ and $Z$ are equal, then any isomorphism between
their respective sets of primitive cycles which respects their colors as described above induces an isomorphism between their cofibrant replacements.
\vskip 2mm

\noindent{\bf Remarks.}
Let $X$ be an undirected graph; $X$ is $n$-regular if and only if for every node $x$ of $X$, the cardinal of $X(x,*)$ is $n$.
Remark that an $n$-colored graph is an $n$-regular but not every $n$-regular graph is $n$-colored as shows the snark graph. The Quillen
model that we have just defined in the category of $n$-regular colored graphs cannot be naively extended to the category whose
objects are $n$-regular graphs and the morphisms are the coverings morphisms, since  pushouts do not exist in this category.

\vskip 2mm

We are going to relate $n$-colored graphs to Cayley graphs. Let $G$ be a group and $S$ a set of generators of $G$,
 for every $G$-set $X$, the Cayley graph $C(X,S,G)$ is the directed graph
whose set of nodes is $X$, and for every elements $x$ and $y$ of $X$, the set of arcs between $x$ and $y$ is in bijection with $\{s\in S, s(x) = y\}$.

Let $G_n$ be the group generated by $S_n=\{a_0,...,a_n\}$ such that $a_i^2 =1$ for every $i=0,...,n$. For every $G_n$-set $X$,
the Cayley graph $C(X,S_n,G_n)$ is endowed with the structure of an undirected regular colored graph defined as follows: let $x$ be a node of $X$, there exists an half edge between $x$ and $a_i(x)$ colored by $a_i$. The symmetric of this half edge is the half edge defined by $a_i(a_i(x))=x$.
 We denote by $UC_{G_n}$  the category of $G_n$-sets.

\medskip

\noindent{\bf Proposition 5.6.}
{\it
The correspondence $C(S_n,G_n)$
which associates to $X$, $C(X,S_n,G_n)$ induces an isomorphism between $UC_{G_n}$ and the category of $n+1$-regular colored graphs $C_{B_n}$.}

\vskip 2mm

{\bf Proof.}
We have only to construct the inverse of $C(S_n,G_n)$. Let $X$ be an $n+1$-regular colored graph. We assume that the colors
are labeled by $a_0,..,a_n$. We associate to $X$ its set of nodes   $X_0$ endowed with the action of $G_n$ defined as follows: if $x\in X_0$ and there exists an arc between $x$ and $y$ colored
by $a_i$, we set $a_i(x)=y$.
\vskip 2mm

\noindent{\bf Remarks.}
The  closed model defined on $C_{B_n}$ induces a closed model on $UC_{G_n}$.

Let $LUC_{G_n}$ be the full subcategory of $UC_{G_n}$
such that for every object $X$ of $LUC_{G_n}$, every $x\in X$, and for every $i=0,...,n$, $a_i(x)\neq x$. The
functor $C(S_n,G_n)$ establishes an isomorphism between $LUC_{G_n}$ and the category of $n+1$-regular colored graphs without a loop.

\vskip 2mm

 Closed models can also be defined in others interesting comma categories associated to $UGph$, here is an example:

\vskip 2mm

\noindent{\bf Definition 5.7.}
The category of bipartite $BUGph$ graphs is the comma category $UGph/A_U$. Thus a bipartite graph is a morphism
$f:X\rightarrow A_U$. A morphism between the objects  $f:X\rightarrow A_U$ and $g:Y\rightarrow A_U$ of $UGph_n/A_U$
 is a morphism $h:X\rightarrow Y$ such that $f = g\circ h$.

Consider the undirected graph $D_n$ which has two nodes $0$ and $1$, and such that there exist $n+1$-arcs $a_0,...,a_n$
between $0$ and $1$. The category $BC_n=C_{D_n}$ of coverings of $D_n$ is the category of graphs which are bipartite and $n+1$-colored.
We deduce the existence of a closed model on $BC_n$ obtained by counting objects of $BC_n$ obtained by attaching a forest
to a cycle.

\vskip 6mm
{\bf 6. Closed models by counting and topology.}
\vskip 6mm

We are going to use the closed model defined on $BC_n$ to study the Galoisian
complexes introduced by Ladegaillerie [15]:

\medskip

\noindent{\bf Definitions 6.1.}
Let $S_n^+$ be the oriented standard affine $n$-simplex whose vertices are labeled $A_0,...,A_n$. We denote by $S_n^-$ the corresponding simplex
with the opposite orientation. Let $I$ be a set (not necessarily numerable), and $(S^+_i)_{i\in I}$ a set of examples of $S_n^+$ and
$(S_i^-)_{i\in I}$ the corresponding set of examples of $S_n^-$. The elements of $(S_i^+)_{i\in I}$ are called the direct simplexes, and the elements of $(S_i^-)_{i\in I}$ are called the undirect simplexes. A Galoisian $n$-complex $C$ is obtained
by gluing elements of $(S_i^+)_{i\in I}$ with elements of $(S_i^-)_{i\in I}$ such that the gluing respect the labeling, affine structures
and inverse orientations. Moreover,we suppose that each face of $C$ belongs to exactly two simplexes; one direct and the
other undirect.

A morphism $f:X\rightarrow Y$ between two Galoisian $n$-complexes is a continuous map which sends a direct simplex to a direct simplex,
an undirect simplex to an undirect simplex, respects the labelings, and the affine structures. This defines the category $CG_n$, whose objects are Galoisian
$n$-complexes and the morphisms are morphisms between Galoisian $n$-complexes.

\vskip 2mm

Let $X$ be a Galoisian $n$-complex. We denote by $\Omega_n^+(X)$ the union of elements of $(S_i^+)_{i\in I}$,  by
$\Omega^-_n(X)$ the union of elements of $(S_i^-)_{i\in I}$, and by $\Omega_n(X)$ the union of $\Omega^+_n(X)$ and $\Omega^-_n(X)$. We can define
$s_j,j=0,...,n$ the involution of $\Omega(X)$ such that for an element $S_i^+$ of $\Omega^+_n(X)$,
$a_j(S_i^+)$ is the unique undirected simplex of $\Omega^-_n(X)$ whose $j$-face is identified
with the $j$-face of $S_i^+$. If $G_n$ is the group generated by $\{a_0,...,a_n\}$ with the relations $a_j^2=1, j=0,...,n$,
the proof of Ladegaillerie [15] 1.2 shows that the correspondence $\Omega_n$ between $CG_n$ and the category of $G_n$-sets, which sends a Galoisian  $n$-complex $X$ to the $G_n$-set $\Omega_n(X)$
endowed with the action that we have just defined is an isomorphism between $CG_n$ and the category of $G_n$-sets. Remark
that  the composition of $C(S_n,G_n)\circ \Omega_n$  defines an isomorphism between $CG_n$ and the category of undirected $n+1$-colored bipartite graphs $BC_n$.
This shows the existence of a closed model on $CG_n$. We will denote by $L_n(X)$ the Cayley graph of the $G_n$-set $\Omega_n(X)$,
defined by the generators of $G_n$,  $a_0,...,a_n$.

\medskip

\noindent{\bf Remarks.}
Let $X$ and $Y$ be finite Galoisian complexes, since $L_n(X)$ is a bipartite graph, it does not have loops, we deduce that if there exists a weak equivalence $L_n(X)\rightarrow L_n(Y)$,
 then $L_n(X)$ and $L_n(Y)$ have the same Ihara Zeta function.

\medskip

\noindent{\bf Question.}
Is it possible to provide a geometric interpretation of the coefficient of the Ihara Zeta function of $L_n(X)$ ?

The Galoisian complex $X_0^n$  defined  by two elements ${S_n^+}_0$ and ${S_n^-}_0$ is homeomorphic to the $n$-sphere $S^n$. For every complex
$X$ defined by $\Omega_n^+(X)\bigcup \Omega_n^-(X)$, there exists a morphism of Galoisian complexes  $p:X\rightarrow X_0^n$ which identifies
the elements of $\Omega_n^+(X)$ to ${S_n^+}_0$, and the elements of $\Omega_n^-(X)$ to ${S_n^-}_0$. This map is ramified at a
$(n-2)$-subcomplex of $C$ (see [15] 1.3).

\vskip 2mm

The subgroup of $G_n$ generated by $\{a_ia_n, i=0,..,n-1\}$ is isomorphic to the free subgroup generated by $n$ elements, $F_n$. Let $X$ be a Galoisian $n$-complex
defined by $\Omega(X)=\Omega_n^+(X)\bigcup\Omega_n^-(X)$. The previous action of $G_n$ on $\Omega(X)$ induces an action of $F_n$ on $\Omega_n^+(X)$. The proof of
Ladegaillerie [15] (p.1725-1726) shows this action induces an isomorphism   between $CG_n$ and the category of $F_n$-sets.

 For every $F_n$-set $X$, we can define the Cayley graph $L_n^+(X)$
defined by the set of generators $a_0a_1,...,a_0a_n$. The isomorphism between $CG_n$ and the category
of $F_n$-sets $E_{F_n}$ induces a closed closed model on $E_{F_n}$. Others closed models can be defined on $E_{F_n}$. In the next section we are going to
present a general construction to transfer the closed model of $Gph$, to the category of $G$-sets for any group $G$.

\vskip 2mm
{\bf  Closed model and $G$-sets.}
\vskip 2mm

Let $G$ be a group, consider the category $C_G$ which has only one object that we denote by $*$, we suppose that
 $Hom_{C_G}(*,*)=G$. An object of the category $\hat {C_G}$, of presheaves over $C_G$ is a set $E$, endowed with an action of $G$.  Thus, $\hat{C_G}$
 is the category of $G$-sets.  We are going to transfer  the closed model defined at the section 4 in the category of directed graphs, to the category of $G$-sets.
On this purpose, we  firstly recall the definition of some canonical functors.

\vskip 2mm

The category $C_D$ (see section 4) can be embedded in $\hat{C_D}$ by using the Yoneda embedding as follows: to $0$, we associate the presheaf $\hat 0$ defined by
$\hat 0(0)=Hom_{C_D}(0,0), \hat 0(1) = Hom_{C_D}(1,0)$.We see that $\hat 0$ is $D_D$ the dot graph. To $1$, we associate
the presheaf $\hat 1$ defined by $\hat 1(0) = Hom_{C_D}(0,1),\hat 1(1)=Hom_{C_D }(1,1)$. We remark that $\hat 1$
is the arc graph $A_D$. Let $X$ be a directed graph, we denote by $C_D/X$ the category whose objects are morphisms of presheaves between the objects
of $C_D$ and $X$. The objects of $C_D/X$ are morphisms $D_D\rightarrow X$ and $A_D\rightarrow X$. Thus, the class of objects of $C_D/X$ can be
identified with the union of the set of nodes of $X$, and its set of arcs. Let $f:U\rightarrow X$ and $g:V\rightarrow X$
be two objects of $C_D/X$, a morphism between $f$ and $g$ is a morphism $h:U\rightarrow V$ such that $f = g\circ h$.
Let $a$ be an arc of $X$, there exists a morphism $f_a:A_D\rightarrow X$ whose image is $a$. We also have morphisms $s_a:D_D\rightarrow s(a)\rightarrow X$
and $t_a:D_D\rightarrow t(a)\rightarrow X$, where $s(a)$ and $t(a)$ are respectively the source and the target of $a$. The source and the target
morphisms $D_D\rightarrow A_D$ induces morphisms of $C_D/X$ between $s_a$ and $f_a$ and between $t_a$ and $f_a$. Remark that if
$a$ is a loop these two morphisms are distinct.

Let $A$ be a set of generators of $G$. For every $G$-set $S$, recall that we have denoted by $C(A,G)(S)$ the Cayley graph of $S$ associated to $G$ and $A$.

Let $U$ be the terminal object of the category of $G$-sets. $U$ is the $G$-set which has a unique element $n$. We denote by $B_A$ the Cayley graph of $U$ defined by $A$.
The objects of the category $C_D/B_A$ are the unique morphism $i_a:D_D\rightarrow B_A$ and the morphisms $c_a:A_D\rightarrow B_A$
which sends $A_D$ to the loop of $B_A$ corresponding to $a$. The morphisms of $C_D/B_A$ are the isomorphisms and
the morphisms $\hat s_a:i_a\rightarrow c_a$ induced by $\hat s:\hat 0\rightarrow \hat 1$ and $\hat t_a:i_a\rightarrow c_a$ induced by $\hat t$. We have a functor $F_A:C_D/B_A\rightarrow C_G$
such that $F_A(i_a)=F_A(c_a)= *$, $F_A(\hat s_a) = Id$, $F_A(\hat t_a) = a$.
 In [8] p. 33, Grothendieck defines an equivalence of categories
$e_{B_A}:\hat{C_D/B_A}\rightarrow Gph/B_A$.

\medskip

\noindent{\bf Proposition 6.2.}
{\it The composition: $e_{B_A}\circ \hat {F_A}:\hat {C_G}\rightarrow Gph/B_A$ is the functor which associates to
a $G$-set $S$, the canonical morphism $C(A,G)(S)\rightarrow B_A$ and henceforth $D(A,G)=e_{B_A}\circ \hat {F_A}$,  has a left adjoint.}

\vskip 2mm

\noindent{\bf Proof.}
Let $E^G$ be a $G$-set, we denote by $E$ the image of $E^G$ by the forgetful functor $G$-$Sets\rightarrow Set$. We have:
 $\hat{F_A}(E^G)(i_A) = \hat{F_A}(E^G)(c_a) = E$. The construction in [8] 5.10.1 shows that:
 $$
 e_{B_A}\circ \hat{F_A}(E^G)(0) = E,
 $$
and
 $$
 e_{B_A}\circ \hat{F_A}(E^G)(1)=\bigcup_{a\in A}\hat {F_A}(E^G(c_a)) = \bigcup_{a\in A}E_a,
 $$

where $E_a = E$. This shows that the set of nodes (resp. the set of arcs ) of the Cayley graph of $E^G$
coincide with the set of nodes (resp. the set of arcs) of $e_{B_A}\circ \hat F_A(E)$.

The restriction of $e_{B_A}\circ \hat{F_A}(E^G)( s)$ to $E_a$ is the identity, and the restriction of
 $e_{B_A}\circ \hat F_A(E^G)(t)$ to $E_a$ is the multiplication by $a$.
 This shows that the Cayley graph of $E^G$  is $e_{B_A}\circ \hat F_A(E^G)$.

 We deduce that the functor $D(A,G)$ has a left and right adjoint since $e_{B_A}$ is an equivalence of categories and $\hat {F_A}$
 has left and right adjoint see [8] proposition 5.1.
\vskip 2mm

\medskip

\noindent{\bf Remark.}
 The existence of a left adjoint of $D(A,G)$ can be shown directly, by showing
that $D(A,G)$ commutes with limits. The functor $\hat C_G\rightarrow Gph$ which sends a $G$-set to its Cayley graph
does not have always a left adjoint since it does not  commutes always with limits.

\vskip 2mm

We consider the closed model defined on $Gph/B_A$ obtained by counting the object $c_n\rightarrow B_A$, where $c_n$ is an $n$-cycle.

We present now the transfer theorem that we are going to use see [5] theorem 3.3.
Let $C$ and $D$ be categories in which limits and colimits exist. Suppose that $C$ is endowed with a closed model.
We denote by $W_C$ the class of weak equivalences, $Cof_C$ the class of cofibrations and $Fib_C$ the class of fibrations of this closed model.
Suppose that there exists a functor $F:C\rightarrow D$ which has a right adjoint functor $G$.

We denote $W_D$ the class of morphisms of $D$ such that for every morphism $f\in W_D$, $G(f)$ is a weak equivalence,

We denote $Fib_D$ the class of morphisms of $D$ such that  for every $f\in Fib_D$, $G(f)$ is a fibration,

We denote by $Fib'_D$ the intersection of $W_D$ and $Fib_D$.

An arrow of $D$ is a cofibration if and only if it has the left lifting property with respect to every element of $Fib'_D$. We denote by $Cof_D$ the class of cofibrations.

\medskip

\noindent{\bf Theorem 6.3.}
{\it With the notations above, suppose that $D$ allows the small object argument, and suppose that for every morphism $d$ of $D$ which is a transfinite
 composition of pushouts of coproducts of morphisms $F(c)$ where $c$ is a weak cofibration, $G(d)$ is a weak equivalence. Then, there exists a closed model on $D$ whose class of weak equivalences
is $W_D$, the class of fibrations is $Fib_D$ and the class of cofibrations is $Cof_D$.}

\vskip 2mm

We can deduce the following result:

\medskip

\noindent{\bf Corollary 6.4.}
{\it Let $G$ be a group, $A$ a set of generators of $G$ and $\hat C_G$ the category of $G$-sets. There exists
a closed model on $\hat C_G$ such that the morphism of $G$-sets $f:X\rightarrow Y$ is a weak equivalence if and only if
$D(A,G)(f):D(X,A,G)\rightarrow D(Y,A,G)$ is a weak equivalence.}

\vskip 2mm

\noindent{\bf Proof.}
To show this result, we apply the theorem 6.4 to transfer the model defined on the category of directed graphs at paragraph 4 with the  functor $D(A,G)$. The category of $G$-sets allows the small object argument. Let $F$ be the left adjoint of $D(A,G)$ and  $d$ be a morphism of $\hat C_G$ which is a transfinite composition of
pushouts of coproducts of morphisms $F(c)$, where $c$ is a weak cofibration. Since weak cofibrations in $Gph$ are isomorphisms, we deduce that
$d$ and $D(A,G)(d)$ are isomorphisms.

\vskip 2mm

We are going to apply this construction to the free group $F_n$. Let $A_n=\{a_1,...,a_n\}$ be a set of generators of $F_n$; we denote by $R_n$ the set of $F_n$-sets such that for every $F_n$-set $X$ in $R_n$, $C(A_n,F_n)(X)$ is  obtained by attaching a forest to a sum of cycles.

\medskip

\noindent{\bf Proposition 6.5.}
{\it A  $F_n$-set $X$ is cofibrant for the closed model obtained by transferring the closed model of $Gph/B_A$ to the category of $F_n$-sets with $D(A_n,F_n)$ if and only if it is an element of $R_n$.}

\medskip

\noindent{\bf Proof.}
Let $X$ be a cofibrant $F_n$-set. The source of the cofibrant replacement $X'$ of $D(A_n,F_n)(X)$ is the sum $\sum_{i\in I}X'_i$ of cycles. Let $h_i$ be the restriction of the cofibrant morphism $h:X'\rightarrow D(A_n,F_n)(X)$ to $X'_i$. There exists an element of $R_n$, $Y_i=D(A_n,F_n)(X_i)$, a morphism $g_i:X_i\rightarrow X$ such that $h_i$ is the restriction of $D(A_n,F_n)(g_i)$ to the cycle of $Y_i$. To see this, consider 
 $x^i_{1},...,x^i_{i_n}$ be the nodes of $X'_i$. We suppose that there exists an arc between $x^i_{i_j}$ and
$x^i_{i_{j+1}}$ if $i_j<i_n$ and an arc between $x^i_{i_n}$ and $x^i_1$.  There exists a
generator $a_{i_j}$ such that $a_{i_j}(h(x^i_{i_j})) = h(x^i_{i_{j+1}})$, $i_j<i_n$ and $a_{i_n}(h(x^i_{i_n}) = h(x^i_1)$. The set $F_n$-set $X_i$ is the unique set which contains the elements $x^i_1,..,x^i_{i_n}, a_{i_j}(x^i_{i_j}) = x^i_{i_{j+1}}$, $j<n$ and $a_{i_n}(x^i_{i_n}) = x^i_1$ and the Cayley graph of $X_i$ is obtained by attaching a minimal forest to the cycle $(x^i_1,..,x^i_{i_n})$.  The sum $\sum_{i\in I}g_i$ is a weak equivalence.

 Consider the commutative diagram:

 $$
 \matrix{ \phi &{\longrightarrow} & \sum_{i\in I} X_i\cr\downarrow &&\downarrow \sum_{i\in I}g_i\cr
 X &{\buildrel{Id_X}\over{\longrightarrow}} & X}
 $$

Since $\sum_{i\in I}h_i$ is a weak equivalence or equivalently a weak fibration, we deduce that the existence of a morphism $f:X\rightarrow \sum_{i\in I}X_i$
which fills the previous commutative diagram. This implies that $X$ is an element of $R_n$.

Conversely, let $X$ be an element of $R_n$, let $f:Y\rightarrow Z$ be a weak equivalence or equivalently a weak fibration such that there exists a commutative diagram:
$$
 \matrix{ \phi &{\longrightarrow} & Y\cr\downarrow &&\downarrow f\cr
 X &{\buildrel{g}\over{\longrightarrow}} & Z}
 $$
   Since the source of $D(A_n,F_n)(X)$ is obtained by attaching a forest to a union of cycles, and $f$ is a weak equivalence the morphism $D(A_n,F_n)(g)$ can be lifted
   to a morphism $D(A_n,F_n)(h):D(A_n,F_n)(X)\rightarrow D(A_n,F_n)(Y)$ which is the image by the functor $D(A_n,F_n)$
   of a morphism $h:X\rightarrow Y$ which makes the previous diagram commutes.
\vskip 2mm

\noindent{\bf Remarks.}
   Let $X$ and $Y$ be $F_n$-sets, if $c(X)$ is a cofibrant replacement of $X$, it is also a fibrant replacement of $X$
   since every morphism of $F_n$-sets is a fibration. 
   
   There exists a functor $c:F_n$-$sets\rightarrow F_n$-$sets$ such that $c(X)$
   is a cofibrant replacement of $X$. To construct $c(X)$, consider a cofibrant replacement $c(X)$ of $X$ and suppose that
   every connected component of $C(A_n,F_n)(c(X))$ is not isomorphic to a tree.

\medskip

\noindent{\bf Proposition. 6.6.}
{\it Let $X$ and $Y$ be $F_n$-sets, the set of morphisms $Hom_{Hot}(X,Y)$ between $X$ and $Y$ in the homotopy category is
 $Hom_{F_n-sets}(c(Y),c(Y))$ where $c(Y)$ and $c(Y)$ are respectively cofibrant
replacements of $X$ and $Y$.}

\vskip 2mm

\noindent{\bf Proof.}
We are going to show that the category $Hot_n$ whose objects are $F_n$-sets and such that for every objects
$X$ and $Y$ of $Hot_n$, $Hom_{Hot_n}(X,Y) = Hom_{F_n-sets}(c(X),c(Y))$ is a localization of the category of $F_n$-sets by the class of weak equivalences.
The morphism $f:X\rightarrow Y$ is a weak equivalence between $F_n$-sets if and only if 
$c(f):c(X)\rightarrow c(Y)$ is  a weak equivalence. This is equivalent to saying that $D(F_n,A_n)(c(f))$ is a weak
equivalence. Since the sources of $D(F_n,A_n)(c(X))$ and $D(F_n,A_n)(c(Y))$ are sum of cycles, this is equivalent to say that $D(F_n,A_n)(c(f))$ and $c(f)$ are  isomorphisms.
We deduce that $Hot_n$ is a localization of the closed model that we have defined on the category of $F_n$-sets. We can conclude by using  the remark at the first line of [6] p.29.

\vskip 2mm

\noindent{\bf  Closed model and Dessins d'enfants.}

\vskip 2mm

In this section, we are going to recall the definition of a dessin d'enfant and see how it is a particular case
of the construction above.

 Let $FCG_2$ be the category of finite Galoisian $2$-complexes, and  $X$ an object of $FCG_2$. The morphism $p_X:X\rightarrow S^2$ is a covering of the $2$-sphere ramified at
three elements that we denote $A_0,A_1$ and $A_2$. We can identify $S^2-\{A_0,A_1,A_2\}$ with $C-\{0,1\}$, the complex line  without two points. Let $[0,1]$ be the segment drawn between $0$ and $1$ in the
complex line,  $X_C=p_X^{-1}([0,1])$ is an undirected graph.
Remark that $p$ induces a morphism $X_C\rightarrow [0,1]$; thus $X_C$ is a bipartite graph.

Let $U^0_X=p_X^{-1}(0)$ and $U^1_X=p_X^{-1}(1)$. The fundamental group $F_2$ of   $\pi_1(C-\{0,1\})$ is the free group generated by two elements. We denote by $s_0$ and $s_1$ its generators. Without restricting the generality, we suppose 
that for each $x\in p_X^{-1}(0)$, the monodromy of $s_0$ induces an action on $X_C(x,*)$ and for every $x\in p_X^{-1}(1)$
the monodromy of $s_1$ induces an action on $X_C(x,*)$. This action is nothing but the restriction of the  action of $F_2$ on $\Omega^+(X)$
(see [17]).
A bipartite graph endowed with such an  action of $F_2$ is called
a dessin d'enfant. Conversely, any finite $F_2$-set define a dessin d'enfant. Let $FS_2$ be the category of finite $F_2$-sets; the
 functor $F: FCG_2\rightarrow FS_2$ which associates to a finite Galoisian $2$-complex
the $F_2$-set defined by its dessin d'enfant induces an isomorphism between $FCG_2$ and $FS_2$ (see [17]).

There exists a one to one correspondence between  finite dessin d'enfants and algebraic curves defined over the algebraic closure $\bar Q$ of the field of rational numbers.
The action of the Galois group $Gal(\bar Q/Q)$ on algebraic curves defined over $\bar Q$ induces an action of $Gal(\bar Q/Q)$ on dessins d'enfants.

We denote by $Hot(FS_2)$  the homotopy category of the closed model defined on $FS_2$-sets.

\medskip

\noindent{\bf Remark.}
Let $D_0$ be the Dessin d'enfant whose underlying graph is $A_U$ and such that $F_2$ acts trivially on the nodes. The Cayley graph
$Cal(D_0)$ associated to this action is $B_2$: the graph which has one node $n$ an two loops $a,b$.

We denote
by $D_1$ the dessin d'enfant whose underlying graph is $P_2$. Let $s_1$ and $s_2$ be the generators of $F_2$. We suppose that $s_1$
acts trivially on the arcs of $P_2$ and $s_2$ defines a non trivial involution on them. The Cayley graph $Cal(D_1)$
associated to this action is a directed graph which has two nodes $x$ and $y$, there exists one directed arc $u$  between  $x$ and $y$,
one loop $u_x$ at $x$, one directed arc $v$ between  $y$ and $x$ and one loop $v_y$ at $y$. The morphism
 $f:Cal(D_1)\rightarrow Cal(D_0)$ defined by $f_0(x) = f_0(y) = n$, $f_1(u) = f_1(u_x) =a, f_1(v)=f_1(v_y) =b$ is a weak equivalence
  for the closed model defined on $Gph$ by counting the cycles but not an isomorphism. This implies that $D_0$ and $D_1$ are weak equivalent.

\medskip

\noindent{\bf Questions.}
Is the action of $Gal(\bar Q/Q)$ on the category of dessins d'enfant induces an action of $Gal(\bar Q/Q)$ on  the image of $FS_2$ in the homotopy category
of the closed model defined on the category of $F_2$-sets ?

We have constructed  (see p. 15) a closed model in $CG_2$ induced by the closed model defined on the category of colored $3$-regular graphs. Is the action of $Gal(\bar Q/Q)$ on dessins d'enfant induces an action on the image of the category of finite Galoisian $2$-complexes in the homotopy category of this closed model ?

\vskip 2mm

\noindent{\bf Conflict of Interests}

\noindent The authors  declare  that there is no conflict of interests.\vskip 2mm

{\bf Bibliography.}

\bigskip

[1].T. Bisson and A. Tsemo, A homotopical algebras of graphs related to zeta series, Homology, Homotopy

and its Applications, 10 (2008), 1-13.

[2] Bisson and A. Tsemo, Symbolic dynamic and the category of graphsTheory and Applications of Categories, Vol. 25, No. 22, 2011.

[3] T. Bisson and A. Tsemo, A  Homotopy equivalence of isospectral graphs New York J. Math. 17 (2011) 295{320.

[4] P. Boldi and S. Vigna, Fibrations of graphs, Discrete Math., 243 (2002), 21-66.

[5] S.E. Crans, Quillen closed model structures for sheaves. J. Pure Appl. Algebra 101 (1995), no. 1, 35

[6] W.G. Dwyer and J. Spalinski, Homotopy theories and model categories, in Handbook of Algebraic
Topology, Editor I. M. James, 73-126, Elsevier, 1995.

[7] M. Ferri, Crystallisations of 2-fold branched coverings of S3, Proc. Amer. Math. Soc.
(1979), 271{276.

[8] A. Grothendieck, and al.  "Théorie des topos et cohomologie étale des schémas. Tome 1: Théorie des topos. Séminaire de Géométrie Algébrique du Bois–Marie, 1963–1964 (SGA 4)." Lecture Notes in Mathematics 269 (1972).

[9] Grothendieck, A., Esquisse d'un programme. London Math. Soc. Lecture Note Ser., 242, Geometric Galois actions, 1, 5�, Cambridge Univ. Press, Cambridge, 1997.

[10] Hatcher, Allen, Algebraic topology, Cambridge University Press (2002)

[11] S. Hirschhorn, Model Categories and Their Localizations, Mathematical Surveys
and Monographs, 99, American Mathematical Society, Providence, RI, 2003.

[12]  G. A. Jones and D. Singerman, Theory of maps on orientable surfaces, Proc. London Math.
Soc. 37 (1978), 273{307.

[13] A. Joyal and M. Tierney, Quasi-categories vs Segal spaces, arXiv:math/06/07/820v1

[14] M. Kotani and T. Sunada, Zeta functions of finite graphs, J. Math. Sci. Univ. Tokyo, 7 (2000), 7-25.

[15] Ladegaillerie Y. Complexes Galoisiens. Transactions of the American Mathematical Society. Volume 352, Number 4, Pages 1723{1741

[16] D.G. Quillen, Homotopical Algebra, Lecture Notes in Mathematics no. 43, Springer-Verlag, Berlin, 1967.

[17] L. Schneps: Dessins d'enfants on the Riemann Sphere, in The Grothendieck Theory of Dessins d'Enfants, LMS Lecture Notes 200, Cambridge U. Press, 1994.

[18]Serre, J-P. Arbres, amalgames, SL(2), Ast閞isque, vol. 46, SMF, Paris, 1977

[19] Terras, Audrey (2010). Zeta Functions of Graphs: A Stroll through the Garden. Cambridge Studies in Advanced Mathematics 128. Cambridge University Press.

[20] Tsemo, A. In preparation.

[21] S. Vigna, A guided tour in the topos of graphs. Technical Report 199-97. Universit`a di Milano.
Dipartimento di Scienze dell扞nformazione, 1997.

\end{document}